\documentclass[12pt]{article} 

\usepackage{amsfonts}
\usepackage{amssymb,amsmath}
\usepackage{amsthm}

\usepackage{todonotes}
\usepackage{adjustbox}

\usepackage[upright]{fourier}\usepackage{baskervald} 
\usepackage{dsfont} 

\usepackage{graphicx}
\usepackage{float}
\usepackage{wrapfig}
\usepackage{subfigure}

\usepackage{url}

\usepackage{hyperref}
\usepackage{geometry}

\usepackage{verbatim}

\usepackage[english]{babel}
\usepackage[latin1]{inputenc}

\usepackage{calligra}
\usepackage[T1]{fontenc}

\usepackage{bbm}

\usepackage{mathrsfs}

\usepackage{color}
\usepackage[color,all]{xy}

\usepackage{extarrows}

\usepackage{caption}

\theoremstyle{definition}
\newtheorem{definition}{Definition}[section]

\theoremstyle{plain}

\newtheorem*{conjecture}{Conjecture}

\theoremstyle{remark}

\newtheorem{remark}[definition]{Remark}

\newcommand{\OO}{\mathscr{O}}

\newcommand{\PP}{\mathbb{P}}
\newcommand{\CC}{\mathbb{C}}
\newcommand{\ZZ}{\mathbb{Z}}
\newcommand{\QQ}{\mathbb{Q}}
\newcommand{\FF}{\mathbb{F}}

\newcommand{\Ann}{\text{\upshape{Ann}}}

\newcommand{\codim}{\text{\upshape{codim}}}
\newcommand{\coker}{\text{\upshape{coker}}}

\newcommand{\Sym}{\text{\upshape{Sym}}}
\newcommand{\Cliff}{\text{\upshape{Cliff}}}

\newcommand{\cL}{\mathscr{L}}

\newcommand{\cM}{\mathscr{M}}

\makeatletter
\newcommand{\subjclass}[2][2010]{%
  \let\@oldtitle\@title%
  \gdef\@title{\@oldtitle\footnotetext{#1 \emph{Mathematics subject classification:} #2}}%
}
\newcommand{\keywords}[1]{%
  \let\@@oldtitle\@title%
  \gdef\@title{\@@oldtitle\footnotetext{\emph{Key words and phrases:} #1.}}%
}
\makeatother

\newcommand{\Christian}[1]{}

\AtEndDocument{\bigskip{
  \textsc{Universit\"at des Saarlandes, Campus E2 4, D-66123 Saarbr\"ucken, Germany} \par  
  \textit{E-mail address}, C.~Bopp: \texttt{bopp@math.uni-sb.de} \par
  \textit{E-mail address}, F.-O.~Schreyer: \texttt{schreyer@math.uni-sb.de} \par
}}

\begin{document}

\title{A version of Green's conjecture in positive characteristic}
\author{Christian Bopp and Frank--Olaf Schreyer}
\date{\today}

\keywords{Canonical curves, syzygies,  Green's conjecture, finite fields}
\subjclass{14Q05, 14H51, 13D02}

\maketitle 

\begin{abstract}
Based on computeralgebra experiments we formulate a refined version of Green's conjecture and a conjecture of Schicho-Schreyer-Weimann which conjecturally also holds in positive characteristic. The experiments are done by  using our Macaulay2 package, which constructs random canonically embedded curves of genus $g\leq 15$ over arbitrary small finite fields. 
\end{abstract}


\section{Introduction}

Green's conjecture on canonical curves has been a challenging problem in algebraic geometry for over $30$ years. It relates the existence of certain pencils of divisors on an algebraic curve defined over $\CC$ to the non-vanishing of certain Betti numbers in minimal free resolution of the canonical model $C\subset \PP^{g-1}$. 
We recall the precise statement. 
The \emph{Clifford-index} of a line bundle $\cL=\OO(D)$ (resp. a divisor $D$) on a smooth projective curve $C$  of genus $g$ is defined as
$$
\Cliff (\cL)=\deg \cL-2(h^0(C,\cL)-1)=g+1-h^0(C,\cL)-h^1(C,\cL),
$$
 and the Clifford index of a curve $C$ is defined by taking the minimum of all "relevant" line bundles: 
$
\Cliff (C)=\min \big\{\Cliff (\cL) \ \big | \ h^0(C,\cL)\geq 2 \text{ and }h^1(C,\cL)\geq 2  \big\}.
$

If $C$ is non-hyperelliptic of genus $g$, then $C$ can be embedded in the projective space $\PP^{g-1}$ via the map 
$$
\phi_{\omega_C}: C\hookrightarrow C\subset \PP^{g-1},
$$
induced by the canonical bundle $\omega_C$ on $C$.

We denote by $S_C=S/I_C$ the coordinate ring of the canonically embedded curve where $S=\Sym H^0(C,\omega_C)$ is the polynomial ring in $g$ variables. The curve $C\subset \PP^{g-1}$ or more precisely its coordinate ring $S_C$ has a minimal free resolution of length $g-2$,
$$
0 \leftarrow S_C \leftarrow S\ \leftarrow \bigoplus_j S(-j)^{\beta_{1,j}} \leftarrow \bigoplus_j S(-j)^{\beta_{2,j}}
\leftarrow \cdots  \leftarrow \bigoplus_j S(-j)^{\beta_{g-2,j}},
$$
and the $\beta_{i,j}=\beta_{i,j}(C)$ appearing in the resolution above are called the graded Betti numbers of the canonically embedded curve $C$. 
The coordinate ring $S_C$ is Gorenstein and hence the resolution is self-dual and the corresponding Betti table symmetric (i.e., $\beta_{i,i+1}(C)=\beta_{g-2-i,g-i}(C)$).
The  Betti table $(\beta_{i,i+j}(C))_{i,j}$ has the form
	\begin{center}
	\begin{tabular}{c | c c c c c c c}
		& $0$ & $1$ & $2$ & $\dots$ & $g-4$  & $g-3$ & $g-2$ \\
		\hline
		$0$ & $1$ & - & -& $\dots$ & -  &- &- \\
		$1$ & - & $\beta_{1,2}$ &$\beta_{2,3}$ &$\dots$ &  $\beta_{g-4,g-3}$ & $\beta_{g-3,g-2}$ & - \\
		$2$ & - & $\beta_{1,3}$ & $\beta_{2,4}$ & $\dots$ &  $\beta_{g-4,g-2}$ & $\beta_{g-3,g-1}$ & - \\
		$3$ & - & - & -& $\dots$ & -  &- &$1$ \\	
	\end{tabular}.
\end{center}

In \cite{G}, Mark Green conjectured the following relation between the vanishing Betti numbers and the existence of certain special linear series on $C$:
\begin{conjecture}[Green's Conjecture]
	Let $C\subset \PP^{g-1}$ be a smooth canonically embedded curve of genus $g$ defined over a field of characteristic $0$. Then 
	$$
	\beta_{i,i+2}(C)=0 \ \text{for all } i\leq p \Longleftrightarrow \Cliff(C)>p.
	$$
\end{conjecture}

	The direction "$\Rightarrow$" was proved  by Green and Lazarsfeld in the appendix to \cite{G}, and the other direction was proved for general curves in Voisin's landmark papers \cite{V1} and \cite{V2}. 
	Recall that the Clifford index of a general curve is known to be $m=\big \lceil \frac{g-2}{2} \big \rceil$. Thus, in this case Green's conjecture predicts that $\beta_{m-1,m+1}(C)=0$. 
	
	By now, several other cases of Green's conjecture have been established, although the conjecture is still open in full generality. For instance, Aprodu showed in \cite{Aprodu05} that Green's conjecture holds for general $k$-gonal curves 
	(i.e., curves admitting admitting a one dimensional linear series of minimal degree $k$)
	and more recently Aprodu and Farkas  showed in \cite{AF11} that Green's conjecture holds for every smooth curve on an arbitrary $K3$ surface.

   However, for curves defined over a field of positive characteristic, it is known that Green's conjecture fails in some cases. By \cite{Sch} and \cite{Muk10},  the Betti tables of a general curve of genus $7$ (resp. 9) defined over a field of characteristic $2$ (resp. 3) have the following Betti tables
   \begin{center}	
   	{\small{
   			\begin{tabular}{|cccccc}
   				\multicolumn{1}{@{}}{\ }&& \\[-2mm] \hline
   				1 & . & .& .& .& . \\
   				.& 10&16 &1 &.& .\\
   				.& .&1&16& 10 & .  \\
   				.&  .&.& .&.& 1				
   				\\[-2mm]
   				\multicolumn{1}{@{}}{\ }&&
   			\end{tabular}
   	}}
   	\ \ 	and		\ \ 
   	{\small{
   			\begin{tabular}{|cccccccc}
   				\multicolumn{1}{@{}}{\ }&& \\[-2mm] \hline
   				1 & . & .& .& .& .& .& . \\
   				.& 21&64 &70 &6& .& .& .\\
   				.& .&.&6& 70 & 64 & 21 & . \\
   				.&  .&.& .&.& .& .&1				
   				\\[-2mm]
   				\multicolumn{1}{@{}}{\ }&&
   			\end{tabular}
   	}}, respectively.
   \end{center}
As usual we denote by $g^r_k\subset |\OO(D)|$ an $r$-dimensional linear series of divisors of degree $k$.    
   We further denote by $\cM^1_{g,k}\subset \cM_g$  the gonality stratum inside the moduli space of curves of genus $g$, consisting (set-theoretically) of curves having a one dimensional linear series of degree $k$ (i.e., a $g^1_k$). 
   
   A canonically embedded curve $C\subset \PP^{g-1}$ which represents a general point in $\cM^1_{g,k}$ (for $k<\big\lceil \frac{g}{2}\big \rceil$) has Clifford index $(k-2)$ and lies on a $(k-1)$-dimensional rational normal scroll $X$ of degree $g-k+1$. This scroll $X$ is swept out by the unique $g^1_k$ on $C$ and contributes with an Eagon-Northcott complex of length $(g-k)$ (with $\beta_{g-k-1,g-k+1}(X)=g-k$) to the minimal free resolution of the curve $C\subset \PP^{g-1}$ (see e.g., \cite{Sch}). In \cite{SSW} the following conjecture is made.
   
   \begin{conjecture}[Schicho-Schreyer-Weimann]
   	Let $C\subset \PP^{g-1}$ be a smooth canonically embedded curve  of genus $g\neq 6$ and let $k<\big\lceil \frac{g}{2}\big \rceil$. Then $W^1_k(C)$ is a reduced single point if and only if
     $\beta_{k-2,k}(C)=g-k$ and $\beta_{i,i+2}(C)=0$ for $i<k-2$.
   \end{conjecture}
	In case $g=6$, smooth canonical curves with $\beta_{1,3}\neq 0$ have $\beta_{1,3}=3$ and are either trigonal or isomorphic to a  smooth plane quintic, by Petri's theorem.
	
   In characteristic $0$,  the conjecture above has  recently been proven by Farkas and Kemeny \cite{FK16}. However, similarly to Greens conjecture, there are exceptional cases to the conjecture above in positive characteristic. For instance, Sagraloff showed that the conjecture above does not hold in characteristic $3$ for a general $C\in \cM^1_{9,5}$ (see  \cite[Rem. 4.5.3]{Sag}).
   
   Based on computer experiments using \emph{Macaulay2} \cite{M2} we suggest the following refinement of the classical generic Green's conjecture and the conjecture above.
   \begin{conjecture}[Refined Green Conjecture]\label{Conj_RGC}
   	Let $C\subset \PP^{g-1}$ a canonically embedded curve defined over an algebraically closed field  and  let 
   	$$
   	\text{strand}_2(S_C): \ 0 \leftarrow S(-3)^{\beta_{1,3}} \xleftarrow{\varphi_2} S(-4)^{\beta_{2,4}} \xleftarrow{\varphi_3} \dots  \xleftarrow{\varphi_{g-3}} S(-(g-1))^{\beta_{g-3,g-1}} \leftarrow 0
   	$$	
   	be the second linear strand of a minimal free resolution of the coordinate ring $S_C$ (here $S(-(i+2))^{\beta_{i,i+2}}$ is in homological degree $i$).
   	Then
   	\begin{enumerate}
   		\item[(a)] 	$H_i(\text{strand}_2(S_C))$ is a module of finite length for all $i\leq p$ if and only if $ \Cliff(C)>p$. 
   		\item[(b)] If $C$ is general inside the gonality stratum $\cM^1_{g,k}\subset \cM_g$ with $2<k <\big \lceil\frac{g+2}{2} \big \rceil$ then $H_{k-2}(\text{strand}_2(S_C))$ is supported on the rational normal scroll swept out by the unique $g^1_k$ on $C$.
   	\end{enumerate}	
   \end{conjecture}

\begin{remark}\label{rem_GC_RGC}
	\begin{enumerate}
		\item[(i)] 		
		Note that the generality assumption in part (b) is necessary. For instance, the statement will not hold for a curve $C\in \cM^1_{g,k}$ with $n\geq 2$ independent pencils of divisors of degree $k$. In this case we expect the support of $H_{k-2}(\text{strand}_2(S_C))$ to be the union of the $n$ rational scrolls defined by the $n$ pencils.
		\item[(ii)] If Green's conjecture holds for a smooth canonical curve $C\subset \PP^{g-1}$ then also part (a) of the refined Green conjecture holds: 		
		A single Green-Lazarsfeld syzygy 
		$$S^{\beta_{g-3-k,g-2-k}}(-(g-2-k)) \overset{s}{\longleftarrow}S^1(-(g-1-k))$$
		  coming from a linear series $|L|$ computing $\Cliff(C)=k-2$ has order ideal generated by at most 
		  $h^0(C,L)+h^0(C,\omega_C\otimes L^{-1})-1=g-1-\Cliff(C)$
		   forms (see \cite[Appendix]{G}). 
		 The corresponding row 
		 $$S^1(-k) \overset{s^t}{\longleftarrow}S^{\beta_{k-1,k+1}}(-(k-1))$$
		 is the presentation matrix of a quotient of $H_{k-2}(\text{strand}_2(S_C))$.
		 Thus,  the zero locus of the annihilator of $H_{k-2}(\text{strand}_2(S_C))$  contains a linear subspace of dimension at least $\Cliff(C)-1\geq 0$.			
	\end{enumerate}
\end{remark}

	\paragraph*{Acknowledgment} 
We would like to thank Mike Stillman for valuable Macaulay2 advices.	
This work is a contribution	to Project 1.7 of the SFB-TRR 195 "Symbolic Tools in Mathematics and their Application" of the German Research Foundation (DFG).


\section{Testing the refined Green conjecture}
As we have already stated in the introduction, the refined Green conjecture is based on experiments using \emph{Macaulay2}.
All the experiments are done by using our \emph{Macaulay2}-package \cite{M2-RandomCurvesOverVerySmallFiniteFields}. The package computes ideals of random canonically embedded curves of genus $g\leq 15$ over the field $\FF_p=\ZZ/p\ZZ$ for some prime number $p$.  
The construction of such random curves for genus $g\leq 14$ goes along the lines of the unirationality proofs of $\cM_g$ for $g\leq 14$ (see \cite{AC83}, \cite{S81}, \cite{CR} and \cite{V05}). For an implementation of these constructions see  \cite{ST02} and \cite{S11} and \cite{M2-RandomCurves}. For genus $15$ the construction follows \cite{Sch-g15} and the corresponding \emph{Macaulay2}-package \cite{M2-MatFac15}.

The new achievement in our package is that it works over arbitrary small fields. 
Since the unirational parametrization of $\cM_g$ is only a rational map, bad choices of parameters in the construction steps might end up in the indeterminacy locus or other undesired subloci. Roughly speaking, the chances to end up in a bad locus of codimension $c$ in one of the construction steps is about  ${1}:{p^c}$, if we work over $\FF_p$. However, for $\FF_2$ this is only very rough: For example, $90\%$ of all the hypersurfaces defined over $\FF_2$ (inside some projective space) correspond to points in the discriminant locus of singular hypersurfaces (see \cite{GvBS}). Since our constructions proceeds in several steps, each of which has a certain failure rate, finding a single smooth canonical curve over a very small field might need a lot of attempts.
The functions in the package \cite{M2-RandomCurvesOverVerySmallFiniteFields} catch all possible  missteps in the construction and try again until success. 

\medskip

The actual test of the statements in the refined Green conjecture is straight forward.
In the first step we test for each $g\leq 15$ whether a prime number $p\leq 101$ is a possible exceptional characteristic for the generic Green conjecture, i.e., for none of these curves the critical Betti number $\beta_{m-1,m+1}$ (with $m=\big \lceil \frac{g-2}{2} \big \rceil$) vanishes. 
This yields the following conjectural exceptional cases to the generic Green conjecture.

\begin{center}  	
	\begin{tabular}{|c|c|c|}
		\hline
		genus & char$(\mathbb{F}_p)$ & conjectural generic extra syzygies \\ \hline \hline
		$7$ & $2$  & $\beta_{2,4}=1$\\ \hline
		$9$ & $3$ & $\beta_{3,5}=6$\\ \hline
		$11$ & $2$, $3$ & $\beta_{4,6}=28,\ 10$ \\ \hline
		$12$ & $2$ & $\beta_{4,6}=1$ \\ \hline 
		$13$ & $2$, $5$& $\beta_{5,7}=64,\ 120$ \\ \hline
		$15$ & $2$, $3$, $5$ & $\beta_{6,8}= 299,\ 390,\ 315$ \\ \hline			
	\end{tabular}
\captionof{table}{Conjectural exceptional cases for the generic Green conjecture in positive characteristic and $g\leq 15$.}	\label{excCases} 
\end{center}

\begin{remark}
	\begin{enumerate}
\item[(a)]		
By semi-continuity on the Betti numbers, our experiments show that the generic Green conjecture holds for $g\leq 15$ and all $p\leq 101$ which are not listed in the table above. This in turn means that the generic Green conjecture holds over  $\QQ$ (and hence over $\CC$) for $g\leq 15$. However, our experiments do neither show that the generic Green conjecture really fails for the cases listed in Table \ref{excCases} nor that there are no further exceptional primes. This second issue has recently been solved
by  Eisenbud and Schreyer (see \cite{ES18}), who showed that the generic Green conjecture holds for $g\leq 15$ and all primes not listed in the table above. In \cite{ES18}, Eisenbud and Schreyer conjecture furthermore that Green's conjecture holds for a general curve of genus $g$ over a field of characteristic $p>0$ if $p\geq \frac{g-1}{2}$. 

 Our evidence, that the cases listed in Table \ref{excCases} are really exceptional cases for the generic Green conjecture is a heuristic argument: The locus of curves with 
 larger Betti numbers than the general curve in $\cM_g$ of the given characteristic is a proper subscheme of $\cM_g$. Thus it is unlikely (but not impossible) that we always end up in this particular subscheme.

\item[(b)]
Even if we do not know that the generic Green conjecture fails for the cases listed in Table \ref{excCases}, the examples which contributed to this table are counterexamples to the "full" Green conjecture. In all these examples the homology at the critical position is a module of finite length. This does not happen if the extra syzygies are induced by a linear series on the curve (see Remark \ref{rem_GC_RGC} above).

For genus $7$ and $12$ this follows also from the fact that the one extra syzygy is too small to be induced by a linear series. 

\item[(c)]
We mention that general genus $10$ or genus $14$ curves defined over a field of characteristic $2$ do in general not have extra syzygies. However, there exist smooth curves of genus $10$ and $14$ which have precisely one extra syzygy. These examples also violate the statement of the "full" Green conjecture in characteristic $2$. 
\end{enumerate}
\end{remark}

For each exceptional pair $(g,p)$ in Table \ref{excCases} we computed $500$ random examples for $g<15$ and $100$ random examples for $g=15$.
In all examples in which a homology group in the second linear strand is a non-trivial finite length module, this homology group is precisely $H_{m-1}((\text{strand}_2(S_C))$ for $m=\big\lceil \frac{g-2}{2} \big\rceil$ and $\beta_{m-1,m+1}(C)$ is the first non-zero Betti number in $\text{strand}_2(S_C)$.
Hence, in order to test the refined Green conjecture, it was in all examples sufficient to study the first non-zero map $\varphi_n$
$$
0 {\longleftarrow } S(-j)^{\beta_{n-1,n+1}} \overset{\varphi_n} {\longleftarrow }S(-j)^{\beta_{n,n+2}}
$$
 in the second linear strand.

The computation of the critical Betti number and the map $\varphi_n$ are not part of our \emph{Macaulay2}-package \cite{M2-RandomCurvesOverVerySmallFiniteFields}. Supplementary \emph{Macaulay2}-files, as well as the complete data of our experiments can be found here: 
\begin{center}
		\href{https://www.math.uni-sb.de/ag/schreyer/index.php/computeralgebra}{https://www.math.uni-sb.de/ag/schreyer/index.php/computeralgebra}.
\end{center}

In order to be able to compute $\varphi_n$ for large genus, we use the newly implemented fast syzygy algorithms (see also \cite{EMSS}).

For each example we test if $M:=\coker (\varphi_n)$ is a module of finite length, by computing the Hilbertfunction of $M$. Since rational normal scrolls are non-degenerate varieties of minimal degree, i.e., varieties satisfying $\codim(\ \_\ )+1=\deg(\ \_\ )$ (see \cite{EH}), computing degree and dimension of $M$ gives a hint whether $\Ann(M)$ contains a rational normal scroll. If moreover $\Ann(M)$ is arithmetically Cohen-Macaulay and is contained in  the ideal of $C$, it is indeed a rational normal scroll (or $g= 6$ and $C$ is isomorphic to a plane quintic). 

\begin{remark}\label{rem_mult_scrolls}
In our experiments we frequently encounter the case that $\deg(M)$ is an integer multiple of $\codim(M)+1$. 
We think that in these cases the curve $C$ 
has several $g^1_k$'s, possibly counted with multiplicities. But, due to limited memory and time resources, we were not able to
compute the annihilator $\Ann(M)$ and its primary decomposition. Hence, we did not verify our interpretation.	
\end{remark}

In Table \ref{table_g11_p2} below  we give the full data of our experiments for the case $(g,p)=(11,2)$. The column "$\#$" tracks how often a particular case occurs and in the column "RGC" we display the degree and the dimension of $\Ann(M)$.
For the first two rows of Table \ref{table_g11_p2} we know the validity of the refined Green conjecture for all tested curves. For curves collected in the other rows, the validity is plausible by Remark \ref{rem_mult_scrolls}  but not proven. 

The complete data for other pairs $(g,p)$ in Table \ref{excCases} looks similar to the $(g,p)=(11,2)$ case. For a complete list of such tables we refer to \cite[Section 6.4]{Bopp-thesis}.

 \begin{remark}
	\begin{enumerate}
		\item[(a)] Computing several examples of odd genus $g$ curves over $\FF_p$, one expects to get a curve inside $\cM^1_{g,k}$ for $k=\big \lceil \frac{g}{2} \big \rceil$ with a chance of roughly $\frac{1}{p}$ (see \cite{GvBS})
		\item[(b)] Note that the last Betti table in Table \ref{table_g11_p2} with critical Betti number $\beta_{4,6}=50$ does not seem to fit in the schematics at first sight. Following the schematics for $\beta_{4,6}$, one might expect that this Betti number is $\beta_{4,6}=48$. But on the other hand, we expect that such curves have precisely $10$  $g^1_6$'s, since $\text{supp}(M)$ is $5$-dimensional and has degree $10\cdot 6=60$ (c.f Remark \ref{rem_mult_scrolls}). 
		
		The experiments in \cite{Sch-pos_char} suggest that in characteristic $0$ or large characteristic,
			a general genus $11$ curve with a $g^2_8$ (whose plane model has the expected $10$ ordinary double points and therefore 	
		$10$ $g^1_6$'s) already has Betti number $\beta_{4,6}=50=5\cdot 10$. This fits with the case over $\FF_2$ above.
	\end{enumerate}
\end{remark}  

\begin{table}[H] \label{table_g11_p2}
	\begin{center}
	\begin{adjustbox}{width=\textwidth}			
		\begin{tabular}{|c|c|c|c|c|}\hline
			genus & char$(\FF_p)$ & \#   & RGC & Betti table \\ \hline
			&  & 230&  $(60,0)$ &
			{\scriptsize{
					\begin{tabular}{|cccccccccc}
						\multicolumn{1}{@{}}{\ }&& \\[-2mm] \hline
						1 & . & .& .& .& .& .& .&.&. \\
						.& 36&160 &315 &288& 28& .& .&.&.\\
						.& .&.&.& 28 & 288 & 315 & 160& 36 & . \\
						.&.&.&  .&.& .&.& .& .&1				
						\\[-2mm]
						\multicolumn{1}{@{}}{\ }&&
					\end{tabular}
			}} \\ \cline{3-5}
			&  &  76 &  $(6,5)$ &
			{\scriptsize{
					\begin{tabular}{|cccccccccc}
						\multicolumn{1}{@{}}{\ }&& \\[-2mm] \hline
						1 & . & .& .& .& .& .& .&.&. \\
						.& 36&160 &315 &288& 30& .& .&.&.\\
						.& .&.&.& 30 & 288 & 315 & 160& 36 & . \\
						.&.&.&  .&.& .&.& .& .&1				
						\\[-2mm]
						\multicolumn{1}{@{}}{\ }&&
					\end{tabular}
			}} \\ \cline{3-5}
			& &  82 & $(12,5)$& 
			{\scriptsize{
					\begin{tabular}{|cccccccccc}
						\multicolumn{1}{@{}}{\ }&& \\[-2mm] \hline
						1 & . & .& .& .& .& .& .&.&. \\
						.& 36&160 &315 &288& 32& .& .&.&.\\
						.& .&.&.& 32 & 288 & 315 & 160& 36 & . \\
						.&.&.&  .&.& .&.& .& .&1				
						\\[-2mm]
						\multicolumn{1}{@{}}{\ }&&
					\end{tabular}
			}} \\ \cline{3-5}	
			& &  55 & $(18,5)$& 
			{\scriptsize{
					\begin{tabular}{|cccccccccc}
						\multicolumn{1}{@{}}{\ }&& \\[-2mm] \hline
						1 & . & .& .& .& .& .& .&.&. \\
						.& 36&160 &315 &288& 34& .& .&.&.\\
						.& .&.&.& 34 & 288 & 315 & 160& 36 & . \\
						.&.&.&  .&.& .&.& .& .&1				
						\\[-2mm]
						\multicolumn{1}{@{}}{\ }&&
					\end{tabular}
			}} \\ \cline{3-5}
			$g=11$&$p=2$ &  24 &  $(24,5)$&
			{\scriptsize{
					\begin{tabular}{|cccccccccc}
						\multicolumn{1}{@{}}{\ }&& \\[-2mm] \hline
						1 & . & .& .& .& .& .& .&.&. \\
						.& 36&160 &315 &288& 36& .& .&.&.\\
						.& .&.&.& 36 & 288 & 315 & 160& 36 & . \\
						.&.&.&  .&.& .&.& .& .&1				
						\\[-2mm]
						\multicolumn{1}{@{}}{\ }&&
					\end{tabular}
			}} \\ \cline{3-5}
			& &  10 & $(30,5)$&
			{\scriptsize{
					\begin{tabular}{|cccccccccc}
						\multicolumn{1}{@{}}{\ }&& \\[-2mm] \hline
						1 & . & .& .& .& .& .& .&.&. \\
						.& 36&160 &315 &288& 38& .& .&.&.\\
						.& .&.&.& 38 & 288 & 315 & 160& 36 & . \\
						.&.&.&  .&.& .&.& .& .&1				
						\\[-2mm]
						\multicolumn{1}{@{}}{\ }&&
					\end{tabular}
			}} \\ \cline{3-5}
			& &  14 &  $(36,5)$&
			{\scriptsize{
					\begin{tabular}{|cccccccccc}
						\multicolumn{1}{@{}}{\ }&& \\[-2mm] \hline
						1 & . & .& .& .& .& .& .&.&. \\
						.& 36&160 &315 &288& 40& .& .&.&.\\
						.& .&.&.& 40 & 288 & 315 & 160& 36 & . \\
						.&.&.&  .&.& .&.& .& .&1				
						\\[-2mm]
						\multicolumn{1}{@{}}{\ }&&
					\end{tabular}
			}} \\ \cline{3-5}
			& &  6 &  $(42,5)$&
			{\scriptsize{
					\begin{tabular}{|cccccccccc}
						\multicolumn{1}{@{}}{\ }&& \\[-2mm] \hline
						1 & . & .& .& .& .& .& .&.&. \\
						.& 36&160 &315 &288& 42& .& .&.&.\\
						.& .&.&.& 42 & 288 & 315 & 160& 36 & . \\
						.&.&.&  .&.& .&.& .& .&1				
						\\[-2mm]
						\multicolumn{1}{@{}}{\ }&&
					\end{tabular}
			}} \\ \cline{3-5}
			& &  2 &  $(48,5)$&
			{\scriptsize{
					\begin{tabular}{|cccccccccc}
						\multicolumn{1}{@{}}{\ }&& \\[-2mm] \hline
						1 & . & .& .& .& .& .& .&.&. \\
						.& 36&160 &315 &288& 44& .& .&.&.\\
						.& .&.&.& 44 & 288 & 315 & 160& 36 & . \\
						.&.&.&  .&.& .&.& .& .&1				
						\\[-2mm]
						\multicolumn{1}{@{}}{\ }&&
					\end{tabular}
			}} \\ \cline{3-5}
			& &  1 & $(60,5)$&
			{\scriptsize{
					\begin{tabular}{|cccccccccc}
						\multicolumn{1}{@{}}{\ }&& \\[-2mm] \hline
						1 & . & .& .& .& .& .& .&.&. \\
						.& 36&160 &315 &288& 50& .& .&.&.\\
						.& .&.&.& 50 & 288 & 315 & 160& 36 & . \\
						.&.&.&  .&.& .&.& .& .&1				
						\\[-2mm]
						\multicolumn{1}{@{}}{\ }&&
					\end{tabular}
			}} \\ \hline																						
		\end{tabular}
		\end{adjustbox}	
	\end{center}
	\begin{footnotesize}\caption{Betti tables of $500$ random examples of genus $11$ curves over $\mathbb{F}_2$}\label{table_g11_p2}\end{footnotesize}
\end{table}


\end{document}